\documentclass[11pt]{article}
\usepackage{amsmath}
\usepackage{amsfonts}
\usepackage{graphicx}
\usepackage{color}
\usepackage[utf8]{inputenc}
\usepackage[T1]{fontenc}
\usepackage[french]{babel}

\makeatletter
\def\@sect#1#2#3#4#5#6[#7]#8{%
  \ifnum #2>\c@secnumdepth
    \let\@svsec\@empty
  \else
    \refstepcounter{#1}%
    \protected@edef\@svsec{\@seccntformat{#1}\relax}%
  \fi
  \@tempskipa #5\relax
  \ifdim \@tempskipa>\z@
    \begingroup
      #6{%
        \@hangfrom{\hskip #3\relax\@svsec}%
          \interlinepenalty \@M #8\@@par}%
    \endgroup
    \csname #1mark\endcsname{#7}%
    \addcontentsline{toc}{#1}{%
      \ifnum #2>\c@secnumdepth \else
        \protect\numberline{\csname the#1\endcsname.}%
      \fi
      #7}%
  \else
    \def\@svsechd{%
      #6{\hskip #3\relax
      \@svsec #8}%
      \csname #1mark\endcsname{#7}%
      \addcontentsline{toc}{#1}{%
        \ifnum #2>\c@secnumdepth \else
          \protect\numberline{\csname the#1\endcsname.}%
        \fi
        #7}}%
  \fi
  \@xsect{#5}}
\def\@seccntformat#1{\csname the#1\endcsname.\quad}
\makeatother

\newtheorem{theo}[equation]{Th\'eor\`eme}
\newtheorem{lem}[equation]{Lemme}
\newtheorem{fait}[equation]{Fait}
\newtheorem{proposition}[equation]{Proposition}
\newtheorem{cor}[equation]{Corollaire}

\makeatletter

\@addtoreset{equation}{section}
\makeatother

\newcommand{\carrenoir}{\rule{0.5em}{0.5em}}
\makeatletter
\newenvironment{demo}[1][\@empty]{\textbf{D\'emonstration~%
\ifx\@empty#1:\else #1~:\fi~}}
{\hfill\carrenoir\nolinebreak\vspace{2mm}}
\makeatother

\newcommand{\oper}[2]{\newcommand{#1}{\mathop{\mathrm{#2}}\nolimits} }
\oper{\Vol}{Vol}
\newcommand{\de}{\mathrm{ d }}
\oper{\Ker}{Ker}
\oper{\Ima}{Im}
\oper{\dimension}{dim}

\oper{\injrad}{inj}
\newcommand{\Z}{\mathbb Z}
\newcommand{\R}{\mathbb R}

\newcommand{\smat}[1]{\left(\begin{smallmatrix}#1\end{smallmatrix}\right)}

\DeclareSymbolFont{greek}{OML}{ptmcm}{m}{it}
\DeclareMathSymbol{\codiff}{\mathord}{greek}{"0E}
\DeclareMathSymbol{\prodint}{\mathord}{greek}{"13}

\title{Minoration du spectre des variétés hyperboliques de dimension~3}
\author{Pierre Jammes}
\date{}
\begin{document}

\maketitle
{\small 
\textsc{Résumé.---}
Soit $M$ une variété hyperbolique compacte de dimension~3, de diamètre~$d$
et de volume $\leq V$. Si on note $\mu_i(M)$ la $i$-ième valeur propre
du laplacien de Hodge-de~Rham agissant sur les 1-formes coexactes de $M$,
on montre que $\mu_1(M)\geq \frac c{d^3e^{2kd}}$ et 
$\mu_{k+1}(M)\geq \frac c{d^2}$, où $c>0$ est une constante ne dépendant 
que de $V$, et $k$ est le nombre de composantes connexes de la partie
mince de $M$. En outre, on montre que pour toute 3-variété hyperbolique
$M_\infty$ de volume fini avec cusps, il existe une suite $M_i$ de
remplissages compacts de $M_\infty$, de diamètre $d_i\to+\infty$ telle que 
et $\mu_1(M_i)\geq \frac c{d_i^2}$.

Mots-clefs : laplacien de Hodge-de~Rham, formes différentielles, variétés
hyperboliques.

\medskip
\textsc{Abstract.---}
Let $M$ be a compact hyperbolic 3-manifold of diameter $d$ and volume $\leq V$.
If $\mu_i(M)$ denotes the $i$-th eigenvalue of the Hodge laplacian acting
on coexact 1-forms of $M$, we prove that $\mu_1(M)\geq \frac c{d^3e^{2kd}}$
and $\mu_{k+1}(M)\geq \frac c{d^2}$, where $c>0$ depends only on $V$,
and $k$ is the number of connected component of the thin part of $M$.
Moreover, we prove that for any finite volume hyperbolic 3-manifold 
$M_\infty$ with cusps, there is a sequence $M_i$ of compact fillings
of $M_\infty$ of diameter $d_i\to+\infty$ such that 
$\mu_1(M_i)\geq \frac c{d_i^2}$.


}

\section{Introduction}
L'objet de cet article est d'étudier le comportement des premières
valeurs propres du laplacien de Hodge-de~Rham, qui agit sur les formes
différentielles,  sur une suite de variétés hyperboliques compactes qui 
dégénère à volume majoré. Plus précisément on va donner une minoration 
du spectre en fonction du diamètre de la variété. En dimension~2, le spectre 
du laplacien de Hodge-de~Rham se déduit du spectre du laplacien usuel
agissant sur les fonctions, et en dimension supérieure ou égale à~4, à
volume majoré, il n'existe qu'un nombre fini de variétés hyperboliques 
compactes, le problème se réduit donc à la dimension~3.
 
 Rappelons quelques résultats connus sur le spectre des variétés 
hyperboliques, en commençant par le laplacien agissant sur les fonctions.
Pour une suite de métriques hyperboliques sur une surface compacte donnée
convergeant vers une variété avec un ou plusieurs cusps, B.~Colbois 
et G.~Courtois ont montré dans \cite{cc89} qu'en restriction à l'intervalle 
$[0,\frac14[$, les valeurs propres convergent vers le spectre de la surface 
limite. En dimension~3 on ne peut pas déformer la métrique d'une variété
compacte dans la classe des métriques hyperboliques, mais W.~Thurston
a montré qu'une variété avec cusps était limite d'une suite de variétés
hyperboliques compactes (voir les rappels fait en section~\ref{topo}). 
Dans ce contexte, Colbois et Courtois \cite{cc91} montrent qu'on a 
convergence du spectre dans l'intervalle $[0,1[$. Dans \cite{pf05}, 
F.~Pf\"affle montre un résultat analogue pour l'opérateur de Dirac en 
dimension~2 et~3. Par ailleurs, R.~Schoen donne dans \cite{sc82} en toute 
dimension supérieure ou égale à~3 une minoration uniforme du spectre du 
laplacien agissant sur les fonctions à volume majoré (voir aussi~\cite{dr86}, 
\cite{do87} et~\cite{bcd93}).

Dans le cas du laplacien de Hodge agissant sur les $p$-formes, une variété 
hyperbolique non compacte de dimension~$n$ possède un spectre continu
formé de la demi-droite $[\frac{(n-2p-1)^2}4,+\infty[$, $p<\frac n2$.
Si $n=3$, le spectre continu pour les 1-formes est donc $[0,+\infty[$,
on s'attend donc à ce que des valeurs propres tendent vers 0 quand une
suite de variétés compactes tend vers une variété non compacte.  On sait
que c'est effectivement le cas (\cite{cc90}, théorème~0.3), et J.~McGowan
et J.Dodziuk ont étudié plus en détail ce phénomène dans \cite{mc93} et
\cite{dm95}.

Avant d'énoncer les résultats, précisons quelques notations. Sur une
variété riemannenne $(M,g)$ compacte sans bord, le laplacien
de Hodge-de~Rham agissant sur les $p$-formes est défini par $\Delta=
\de\codiff+\codiff\de$ où $\de$ désigne la différentielle extérieure et
$\codiff$ la codifférentielle, adjoint $L^2$ de $\de$. C'est un opérateur
positif dont le noyau est canoniquement isomorphe à la cohomlogie 
de Rham $H^p(M)$ de la variété. Le cas $p=0$ correspond au laplacien
usuel sur les fonctions, dont on notera $0<\lambda_1(M,g)\leq\lambda_2(M,g)
\leq\ldots$ le spectre non nul. Pour $p=1$, on a la décomposition de 
Hodge $\Omega^1(M)=\de\Omega^0(M)\oplus\Ker\Delta\oplus\codiff\Omega^2(M)$,
le spectre non nul est donc la réunion du spectre $(\lambda_i(M,g))$
du laplacien restreint à $\de\Omega^0(M)$ et du spectre restreint à 
$\codiff\Omega^2(M)$ qu'on notera $0<\mu_1(M,g)\leq\mu_2(M,g)\leq\ldots$.
En dimension~3, par dualité de Hodge, le spectre en degré 2 et 3 est le 
même qu'en degré 1 et 0 respectivement, on est donc ramené à l'étude des 
$(\mu_i(M,g))$.

Avec ces notations, on peut énoncer les résultats de J.~McGowan comme suit:
\begin{theo}[\cite{mc93}]
Il existe des constantes $C_1,C_2,C_3>0$ et $V>0$ telles que si $M,g$ est une
variété hyperbolique compacte de dimension~3 de diamètre $d$ et de volume 
inférieur à $V$, alors
\begin{enumerate}
\item $\#\left\{\displaystyle\mu_i(M,g)\leq\frac1{C_1V(1+d^2)}\right\}
\leq C_2V$;
\item $\#\left\{\displaystyle\mu_i(M,g)\leq\frac1{d^2}\right\}\geq C_3V$.
\end{enumerate}
\end{theo}
Ce résultat a été précisé dans \cite{dm95} en donnant une estimée
du nombre de valeurs propres dans un intervalle $[0,x]$ pour une
suite de variétés qui dégénère en fonction de la géométrie des différentes 
composantes de la partie mince de la variété (voir section~\ref{topo}
pour la définition de cette notion), généralisant des résultats 
analogues pour les fonctions dans l'intervalle $[1,1+x]$ (\cite{cd94}), 
mais en laissant ouvert la question d'une minoration explicite de la 
première valeur propre en fonction du diamètre.  

Le but de l'article est d'y apporter une réponse en démontrant deux 
résultats de minoration du spectre. Le premier donne une minoration
générale qui est exponentielle par rapport au diamètre de la variété pour
la première valeur propre, et quadratique pour la ($k+1$)-ième, où
$k$ est le nombre de parties minces:
\begin{theo}\label{intro:th1}
Pour tout réel $V>0$, il existe une constante $c(V)>0$ telle que
si $M$ est une variété hyperbolique compacte de dimension~3 de volume 
inférieur à $V$, de diamètre $d$ et possédant $k$ parties minces, alors 
$\mu_1(M)\geq \frac c{d^3e^{2kd}}$ et $\mu_{k+1}(M)\geq\frac c{d^2}$.
\end{theo}

 Les estimées de J.~Dodziuk et J.~McGowan pouvaient laisser espérer une 
minoration de la première valeur propre qui soit quadratique par rapport 
au diamètre. Le second résultat de l'article consiste à montrer que pour 
certaines variétés, on a effectivement une telle minoration:
\begin{theo}\label{intro:th2}
Pour toute variété hyperbolique $M_\infty$ non compacte de dimension~3 et de 
volume fini, il existe une constante $c>0$ et une suite $(M_i)$ de 
remplissages compacts de $M_\infty$, de diamètre $d_i\to+\infty$ et telle 
que $\mu_1(M_i)\geq\frac c{d_i^2}$.
\end{theo}
Un remplissage compact de $M_\infty$ est une variété hyperbolique
compacte dont la partie épaisse est topologiquement identique à $M_\infty$
(voir section~\ref{topo}), en particulier le volume des $M_i$ est
uniformément majoré.

Si on essaie de généraliser la démonstration de cette minoration à toutes
les variétés hyperboliques, une difficulté liée à l'interaction 
entre les cohomologies des parties mince et épaisse de la variété apparaît.
Cet aspect topologique du problème n'avait pas été mis en évidence
par les travaux de McGowan et Dodziuk.

 Ces résultats laissent la possibilité qu'il existe pour chaque partie
mince une petite valeur propre à décroissance exponentielle par 
rapport au diamètre, mais cette valeur propre n'apparaîtrait que pour
des topologies particulières. Comme en dimension~2, l'existence d'une
partie mince n'implique pas forcément l'existence d'une valeur
propre exponentiellement petite par rapport au diamètre.

Les résultats de \cite{mc93} et \cite{dm95} reposent sur une technique 
de minoration du spectre (dont l'idée remonte à Cheeger) faisant
intervenir un recouvrement de la variété par des ouverts (\cite{mc93}, 
lemme~2.3). Cette technique échouait à minorer la 
première valeur propre si la cohomologie des intersections des ouverts
du recouvrement était non triviale. Les théorèmes~\ref{intro:th1} 
et~\ref{intro:th2} reposent sur une amélioration du lemme de McGowan 
qui permet de minorer cette première valeur propre. Dans une autre
direction, le principe d'utiliser des recouvrements pour minorer le 
spectre du laplacien de Hodge-de~Rham a déjà conduit à des résultats
très généraux (\cite{ct97}, \cite{ma08}), mais peu précis quand le
rayon d'injectivité est petit, en particulier dans
le contexte des variétés hyperboliques de dimension~3.
 
L'article est organisé comme suit : on fait dans la section~\ref{topo} 
des rappels sur la géométrie et la topologie des variétés hyperboliques
de dimension~3, la section~\ref{lemmes} est consacrée à des minorations
du spectre sur les parties mince et épaisse, la généralisation du
lemme de McGowan est montré dans la section~\ref{mc} et pour finir la
démonstration des théorèmes~\ref{intro:th1} et~\ref{intro:th2} est
donnée dans la section~\ref{minoration}. 

Je remercie B.~Colbois et F.~Naud d'avoir attiré mon attention sur ce problème.

\section{Topologie, géométrie et cohomologie des variétés hyperboliques 
de dimension~3}
\label{topo}
\subsection{Géométrie et topologie des variétés hyperboliques}
Nous allons rappeler ici plusieurs aspects de la topologie et de la
géométrie des variétés hyperboliques de dimension~3 qui interviendront
dans la démonstration des théorèmes \ref{intro:th1} et \ref{intro:th2}.
On se référera principalement à \cite{gr81} et \cite{bp92}.

 Pour toute variété hyperbolique de dimension~3 de volume fini et tout 
réel $a>0$, on notera $M_a=\{x\in M, \injrad(x)>a\}$ l'ensemble des 
points de $M$ dont le rayon d'injectivité est plus grand que $a$.

Selon le lemme de Margulis, il existe une constante $c_{\mathcal M}>0$
telle que la topologie et la géométrie de la partie mince de $M$,
c'est-à-dire le complémentaire $M_m$ de $M_{c_{\mathcal M}}$, soient simples. 
Plus précisément, la partie mince $M_m$ possède un nombre fini de 
composantes connexes, chacune étant soit un cusp (isométrique à 
$T^2\times[0,+\infty[$ muni d'une métrique $g=e^{-2x}g_T
\oplus\de x^2$, où $g_T$ est une métrique plate sur le tore $T^2$),
soit un voisinage tubulaire $\mathcal T$ d'une géodésique fermée 
de longueur $l$. Ce tube $\mathcal T$ est difféomorphe à $B^2\times S^1$
et  isométrique au produit 
$B^2(R)\times[0,l]$ muni de la métrique 
\begin{equation}\label{topo:g}
g=\cosh^2 r\ \de t^2+\de r^2+\sinh^2 r\ \de\theta^2
\end{equation}
où $t$ désigne la coordonnée le long de la géodésique, $r$ la coordonnée
radiale sur la boule $B^2$ et $\theta$ la coordonnée angulaire,
les bord $B\times\{0\}$ et $B\times\{l\}$ étant identifiés par une isométrie
de $B^2$. 

Si on se donne une borne $V$ sur le volume, il n'existe qu'un nombre fini 
de types de difféomorphisme possibles pour la partie épaisse de la variété
(\cite{bp92}, théorème~E.4.8). De plus, si on considère une suite infinie
de variétés compactes distinctes, on peut en extraire une sous-suite 
qui converge vers une variété avec un ou plusieurs cusps.

On peut déduire en particulier de ce qui précède que l'ensemble des 
métriques induites sur une partie épaisse (ou plus généralement
sur $M_a$) de topologie fixée est borné:
\begin{fait}\label{topo:fait}
Soit $a>0$ et $V>0$. Il existe une constante $\tau(V)>0$ telle que si 
$(M,g)$ et $(M',g')$ sont deux variétés hyperboliques de dimension~3 de 
volume inférieur à $V$ et telles que $M_a$ et $M'_a$ soient 
difféomorphes, alors $\frac1\tau g'\leq g\leq \tau g'$ en restriction 
à $M_a$ et $M'_a$.
\end{fait}

 Si on considère une suite de variétés hyperboliques dont la partie épaisse
est de topologie donnée et dont le diamètre tend vers l'infini, la topologie
des éléments de la suite dépend de la manière dont les tubes $\mathcal T
\simeq B^2\times S^1$ sont recollés sur la bord de la partie épaisse.
Pour chaque composante de la partie mince, les recollements possibles
(topologiquement) sont paramétrés par une élément de 
\begin{equation}\label{topo:chir}
\mathcal P=\{(p,q)\in\Z^2,\ p\textrm{ et }q\textrm{ premiers entre eux}\}
\cup\infty
\end{equation}
(voir~\cite{bp92}, section~E.4).
Comme le fait~\ref{topo:fait} s'applique au bord de la partie épaisse,
la géométrie, donc l'aire, du bord du tube reste contrôlée, et converge
quand $R\to+\infty$. Ceci implique que le produit de la circonférence du 
tube avec la longueur tend vers une constante. Compte tenu de l'expression 
de la métrique sur le tube donnée en~(\ref{topo:g}), on a donc
\begin{equation}\label{topo:circ}
le^{2R}\to c^{te} 
\end{equation}
quand $R\to+\infty$.

 Enfin, W.~Thurston a montré (\cite{th97}, théorème~E.5.1 de \cite{bp92}) 
qu'une variété de volume fini possédant $k$ cusps est la limite d'une suite 
de variétés compactes ayant la même partie épaisse (on parle alors de
remplissage compact de la variété non compacte). Plus précisément, comme 
le bord d'un cusp est un tore $T^2$, on peut remplir topologiquement ce cusp 
par un tube $\simeq B^2\times S^1$, chaque remplissage étant paramétré par
un élément $(p_k,q_k)\in\mathcal P$. Thurston montre alors qu'il existe
$K>0$ tel que si $p_k,q_k>K$ pour tout $k$, alors la variété
compacte obtenue peut être munie d'une métrique hyperbolique. En particulier,
si la variété initiale n'a qu'un seul cusp, alors presque tout
remplissage est hyperbolique.

\subsection{Un lemme cohomologique}
 Dans la démonstration des théorèmes~\ref{intro:th1} et \ref{intro:th2},
l'interaction entre les cohomologies des parties mince et épaisse 
de la variété et de leur frontière joue un rôle crucial. Pour comprendre
cette interaction, on aura besoin du lemme suivant qui sera appliqué à
la partie épaisse de la variété.
\begin{lem}\label{topo:image}
Si $M$ est une variété compacte de dimension~3, alors l'image
de l'application $H^1(M)\to H^1(\partial M)$ est de dimension 
$b_1(\partial M)/2$.
\end{lem}
La démonstration s'appuie sur la suite exacte longue
\begin{equation}\label{topo:longue}
\to H^{k-1}(\partial M)\to H_c^k(M)\to H^k(M) \to H^k(\partial M)\to
\end{equation}
dont on peut trouver la démonstration dans~\cite{ta96} (eq.~(9.67) p.~370,
$H_c^k(M)$ désigne la cohomologie à support compact de $M$).
On commence par établir une propriété générale de cette suite :
\begin{lem}
 Les applications $H^k(M) \to H^k(\partial M)$ et $H^{n-k-1}(\partial M)\to 
H_c^{n-k}(M)$ sont les transposées l'une de l'autre pour la dualité de 
Poincaré.
\end{lem}
\begin{demo}
Rappelons d'abord la définition de ces deux applications. La première,
$\Psi_1: H^k(M) \to H^k(\partial M)$, est naturellement induite par 
l'injection $i:\partial M\to M$. La seconde se construit de la manière 
suivante : si $\alpha\in\Omega^{n-k-1}(\partial M)$ est fermée, on se 
donne une forme $\alpha'\in\Omega^{n-k}(M)$ telle que $i^*(\alpha')=\alpha$
et $\de\alpha'=0$ au voisinage de $\partial M$, 
et on pose $\Psi_2([\alpha])=[\de\alpha']$ (on peut vérifier que la
classe de cohomologie de $\de\alpha'$ ne dépend que de celle de $\alpha$,
voir~\cite{ta96}). 

Si $[\alpha]\in H^{n-k-1}(\partial M)$ et $[\beta]\in H^k(M)$, le
calcul du crochet de dualité $\langle\Psi_2([\alpha]),[\beta]\rangle$ donne,
en utilisant la formule de Stokes et le fait que $\beta$ est fermée,
\begin{eqnarray}
\langle\Psi_2([\alpha]),[\beta]\rangle & = & \int_M \de\alpha'\wedge\beta
=\int_{\partial M}i^*(\alpha'\wedge\beta)\nonumber\\
& = & \int_{\partial M}\alpha\wedge i^*(\beta)
=\langle[\alpha],\Psi_1([\beta])\rangle,
\end{eqnarray}
ce qui est précisément l'énoncé du lemme.
\end{demo}

\begin{demo}[du lemme~\ref{topo:image}]
Dans la suite~(\ref{topo:longue}), on considère les deux flèches
\begin{equation}
H^1(M)\to H^1(\partial M)\to H_c^2(M).
\end{equation}
La dimension de l'image de la première est égale à la codimension 
du noyau de la seconde car les deux applications sont transposées 
l'une de l'autre, et elle est aussi égale à la dimension du noyau
de la seconde par exactitude de la suite. Cette dimension est 
donc nécessairement $b_1(\partial M)/2$.
\end{demo}
\section{Spectres des parties minces et épaisse}\label{lemmes}
La minoration du spectre d'une variété hyperbolique se fera en se 
basant sur des minorations des spectres de sa partie épaisse et de sa 
partie mince. Ces minorations seront combinées dans la section suivante 
pour obtenir une minoration globale en utilisant la technique développée 
par McGowan dans~\cite{mc93}.

Dans cette section et la suivante, le spectre des variétés à bord sera
toujours considéré avec les conditions de bord absolues, qui généralise
la condition de Neumann sur les fonctions. Si on note $j:\partial M\to M$
l'injection canonique et $N$ un champ de vecteur normal au bord, 
cette condition peut s'écrire sous l'une des deux formes suivantes :
\begin{equation}
(A)\ \left\{\begin{array}{l}j^*(\prodint_N\omega)=0\\
j^*(\prodint_N\de\omega)=0\end{array}\right.\textrm{ ou }
\left\{\begin{array}{l}j^*(*\omega)=0\\
j^*(*\de\omega)=0\end{array}\right.
\end{equation}
Avec cette condition de bord, le laplacien de Hodge est elliptique 
(cf.~\cite{ta96}) et son spectre est discret. On conservera les notations
$(\lambda_i)$ et $(\mu_i)$ pour désigner ses valeurs propres.

Rappelons d'abord une caractérisation variationnelle du spectre remontant
à J.~Dodziuk et déjà exploitée par J.~McGowan:
\begin{proposition}[\cite{do82}, \cite{mc93}]\label{lemmes:dod}
 Sur une variété compacte sans bord ou avec condition de bord (A), on a
$$\mu_i=\inf_{V_i}\sup_{\varphi\in V_i\backslash\{0\}}\left\{
\frac{\|\varphi\|^2}{\|\omega\|^2},\ \de\omega=\varphi\right\},$$
où $V_i$ parcourt l'ensemble des sous-espaces de dimension $i$ dans
l'espace des $2$-formes exactes lisses.
\end{proposition}
Remarquons que dans cet énoncé, on n'exige pas que les formes $\omega$ et
$\varphi$ vérifient la condition de bord (A).

 J.~Dodziuk en déduit le résultat qui suit et qui nous sera aussi utile:
\begin{cor}[\cite{do82}]\label{lemmes:dod2}
Soit $g$ et $g'$ deux métriques riemanniennes sur une variété compacte $M$
de dimension $n$, et $\tau$ une constante strictement positive. Si les
deux métriques vérifient $\frac1\tau g\leq g'\leq\tau g$, alors
$$\frac1{\tau^{n+3}}\mu_i(M,g)\leq\mu_i(M,g')\leq
\tau^{n+3}\mu_i(M,g)$$
pour tout entier $i>0$.
\end{cor}

Commençons par montrer qu'à volume majoré, le spectre des 1-formes est
uniformément minoré sur la partie épaisse de la variété, et plus généralement
sur $M_a$ pour $a$ donné, ainsi que le spectre des fonctions sur 
$M_a\backslash M_{c_\mathcal M}$:
\begin{lem}\label{lemmes:1}
Pour tout $V>0$ et tout $a\in]0,c_{\mathcal M}[$, il existe des constantes
$c_1(V,a)>0$ et $c_2(V,a)>0$ telles que si $M$ est une variété hyperbolique 
compacte de dimension~3, alors $\mu_1(M_a)\geq c_1$ et 
$\lambda_1(M_a\backslash M_{c_\mathcal M})\geq c_2$.
\end{lem}
\begin{demo}
Si on fixe la topologie de la partie épaisse de la variété, la famille
de métriques à considérer est contrôlée (cf. fait~\ref{topo:fait})
D'après le corollaire~\ref{lemmes:dod2}, la valeur propre 
$\mu_{1,1}(M_a)$ est donc uniformément minorée pour cette famille.
Comme il n'existe qu'un nombre fini de topologies possibles pour la partie
épaisse, la minoration uniforme de $\mu_{1,1}(M_a)$ s'en déduit immédiatement.

Le fait~\ref{topo:fait} reste valable pour $M_a\backslash M_{c_\mathcal M}$, 
donc le raisonnement précédent s'applique aussi à 
$\mu_{0,1}(M_a\backslash M_{c_\mathcal M})$.
\end{demo}

On passe maintenant la la minoration du spectre sur la partie mince.
\begin{lem}\label{lemmes:2}
Il existe une constante $c>0$ telle que si $\mathcal T$ est une composante 
connexe compacte de la partie mince d'une variété hyperbolique de dimension~3,
alors $\mu_{1,1}(\mathcal T)\geq\frac c{R^2}$, où $R$  est le rayon de 
$\mathcal T$.
\end{lem}
La minoration du spectre des parties minces est plus technique. On
utilisera le fait que la métrique sur une composante connexe
de la partie mince est $T^2$-invariante pour se ramener à un problème
unidimensionnel à l'aide du lemme qui suit :
\begin{lem}[\cite{ja04}]\label{lemmes:inv}
Soit $M$ est une variété compacte $M$ dont la métrique $g$ est 
$S^1$-invariante et $L$ la borne supérieure des longueurs des orbites de 
l'action de $S^1$. Si $\lambda$ est une valeur propre du laplacien de 
Hodge-de~Rham vérifiant $\lambda<(2\pi/L)^2$, alors les formes propres
associées sont $S^1$-invariantes.
\end{lem}

\begin{demo}[du lemme~\ref{lemmes:2}]
Rappelons que la métrique sur une composante connexe compacte $\mathcal T$ 
de la partie mince de $M_m$  s'écrit 
\begin{equation}\label{lemmes:g}
g=\cosh^2 r\ \de t^2+\de r^2+\sinh^2 r\ \de\theta^2
\end{equation}
où $t$ désigne la coordonnée le long de la géodésique, $r$ la coordonnée
radiale sur la boule $B^2$ et $\theta$ la coordonnée angulaire.

 L'expression (\ref{lemmes:g}) est indépendante de $t$ et $\theta$,
ce qui signifie que la métrique est invariante sous l'action du
tore $T^2$, paramétré par $t$ et $\theta$. Nous allons d'abord montrer
qu'il suffit d'étudier le spectre en restriction aux formes différentielles
invariantes sous cette action. Pour ce faire, il suffit de décomposer 
l'action du tore en l'action de deux cercles dont la longueur des orbites 
est bornée, pour ensuite appliquer le lemme~\ref{lemmes:inv}.

Les orbites de l'action de $T^2$ sont les courbes de niveau de la fonction
$r$, qu'on notera $T_r$ ($T_r$ est un cercle pour $r=0$ et un tore sinon). 
Pour $r=R$, cette orbite est le bord du tube. Comme sa géométrie est contrôlée
selon le fait~\ref{topo:fait}, son diamètre est majoré en fonction du volume 
de la variété $M$. On peut facilement
décomposer $T^2$ en deux cercles dont les orbites en restriction
à $r=R$ sont bornées par le diamètre du tore, donc par le volume de $M$.
Pour $r<R$, la métrique $g_r$ induite sur $T_r$ par la métrique de $M$
vérifie $g_r<g_R$, donc les longueurs des orbites des actions de $S^1$ 
sont uniformément majorée par celles de $T_R$.

On va minorer le spectre de la partie mince en appliquant le 
lemme~\ref{lemmes:dod}. On peut vérifier que le spectre restreint 
aux formes $T^2$-invariantes peut se calculer avec la formule du lemme
en supposant $\omega$ et $\varphi$ invariantes.

On est donc ammené à étudier le quotient $\mathcal R(\omega)=
\frac{\|\de\omega\|^2}{\|\omega\|^2}$ où $\omega$ est une forme invariante.
En posant $\omega=a_1(r)\de t+a_2(r)\de r+a_3(r)\de\theta$, on a
$\de\omega=a_1'\de r\wedge\de t+a_3'\de r\wedge\de\theta$. 
La continuité de $\omega$ en $r=0$ et le fait qu'elle soit invariante 
imposent que $a_3(0)=0$. Comme on veut minimiser la norme de $\omega$
pour $\de\omega$ fixé, on peut se restreindre au cas où $a_2=0$.
Compte tenu de l'expression (\ref{lemmes:g}) de la métrique, on cherche 
donc à minorer
\begin{equation}
\mathcal R(\omega)=\frac{\displaystyle\int_0^R\left(\frac{\sinh r}{\cosh r}
a_1'^2+\frac{\cosh r}{\sinh r}a_3'^2\right)\de r}
{\displaystyle\inf_{c\in\R}\int_0^R
\left(\frac{\sinh r}{\cosh r}(a_1-c)^2
+\frac{\cosh r}{\sinh r}a_3^2\right)\de r}
\end{equation}

On va maintenant traiter séparément les cas $a_1=0$ et $a_3=0$ avant 
de combiner les deux résultats. On va d'abord étudier le quotient de Rayleigh
\begin{equation}
\frac{\int_0^Ra_1'(r)^2\tanh r\ \de r}
{\inf_{c\in\R}\int_0^R(a_1(r)-c)^2\tanh r\ \de r}.
\end{equation}
Il s'agit d'un quotient de Rayleigh en dimension~1 pour lequel
la mesure de Lebesgue est remplacée par $\tanh(r)\de r$. Malheureusement
la densité de cette mesure dégénère en $0$, ce qui empêche de se 
ramener classiquement à un laplacien de Witten. Ce problème a déjà été
résolu dans \cite{dm95} (section~6) mais on va le traiter ici par des
techniques plus élémentaires. 

On veut éliminer le terme $\tanh r$ dans le quotient de Rayleigh.
Pour préparer cette élimination, et pour se
ramener à un intervalle de longueur fixe, on va commencer par
faire un changement de variable. On pose $I(R)=\int_0^R\tanh^{-\frac23}(r)
\de r$ (cette intégrale est bien convergente et en outre $I(R)\sim R$ quand 
$R\to\infty$), on fait le changement de variable $\de r/\tanh^{\frac23}r
=\de x/I(R)$ et les changements de fonctions $f(x)=a_1(r)$
et $\eta(x)=\tanh^{\frac13} r$, ce qui donne
\begin{eqnarray}
\frac{\int_0^R\left(\frac{\partial a_1}{\partial r}\right)^2\tanh r\ \de r}
{\int_0^R(a_1-c)^2\tanh r\ \de r} & = &
\frac{\int_0^1\left(\frac{\partial f}{\partial x}
\frac{\partial x}{\partial r}\right)^2\eta^5\de x}
{\int_0^1(f-c)^2\eta^5\de x}\nonumber\\
& = & \frac1{I(R)^2}
\frac{\int_0^1f'^2\eta\de x}{\int_0^1(f-c)^2\eta^5\de x}\nonumber\\
& \geq & \frac1{I(R)^2}\frac{\int_0^1f'^2\eta\de x}{\int_0^1(f-c)^2\de x}
\end{eqnarray}
en utilisant pour cette dernière inégalité que $\eta\leq1$. La fonction
$\eta$ dépend de $R$, mais comme $r\geq x(r)$ pour tout $r$ si $R$ est 
suffisamment grand, on a $\eta^3(x)=\tanh r\geq\tanh x$. On est donc
ramené à majorer $\inf_c\int_0^1(f-c)^2\de x$ en fonction de 
$\int_0^1f'^2\tanh^{\frac13}r\ \de x$.

En dimension~1 on a l'inégalité $\inf_c\|f-c\|_\infty\leq \|f'\|_1$ qui
découle de l'inégalité des accroissement finis, et $\emph{a fortiori}$
les inégalités de Sobolev sont vraies quels que soient les exposants.
En particulier on peut écrire, en appliquant une inégalité de Cauchy-Schwarz :
\begin{eqnarray}
\inf_c\|f-c\|_2^2 & \leq & \|f'\|_1^2=\left(\int_0^1|f'(x)|\tanh^{\frac16}x
\cdot\tanh^{-\frac16}x\ \de x\right)^2\nonumber\\
& \leq & \int_0^1f'(x)^2\tanh^{\frac13}x\ \de x\cdot 
\int_0^1\tanh^{-\frac13}x\ \de x
\end{eqnarray}

Finalement, on a bien montré l'existence d'une constante $C_1$ telle que
si $R$ est suffisamment grand, alors
\begin{equation}\label{lemmes:eq1}
\int_0^Ra_1'(r)^2\tanh r\ \de r\geq \frac{C_1}{R^2}
\inf_{c\in\R}\int_0^R(a_1(r)-c)^2\tanh r\ \de r.
\end{equation}

 Le cas de la fonction $a_3$ se traite de manière similaire. Le même 
changement de variable et le changement de fonction $f(x)=a_3(r)$ donne
\begin{eqnarray}
\frac{\int_0^R\left(\frac{\partial a_3}{\partial r}\right)^2\tanh^{-1} r\ \de r}
{\int_0^Ra_3^2\tanh^{-1} r\ \de r} & = & \frac1{I(R)^2}
\frac{\int_0^1f'^2\eta^{-5}\de x}{\int_0^1f^2\eta^{-1}\de x} \geq 
\frac1{I(R)^2}\frac{\int_0^1f'^2\de x}{\int_0^1f^2\eta^{-1}\de x}\nonumber\\ 
& \geq & 
\frac1{I(R)^2}\frac{\int_0^1f'^2\de x}{\int_0^1f^2\tanh^{-\frac13}x\ \de x}
\end{eqnarray}
Comme $f(0)=a_3(0)=0$, on a $\|f\|_\infty\leq\|f'\|_2$, et donc
\begin{equation}
\int_0^1f^2\tanh^{-\frac13}x\ \de x \leq \|f\|_\infty^2
\int_0^1\tanh^{-\frac13}x\ \de x\leq C\cdot\int_0^1f'^2\de x
\end{equation}
On en déduit l'existence d'une constante $C_2$ telle que
si $R$ est suffisamment grand, alors
\begin{equation}\label{lemmes:eq2}
\int_0^Ra_3'(r)^2\tanh^{-1} r\ \de r\geq \frac{C_2}{R^2}
\int_0^Ra_3(r)^2\tanh^{-1} r\ \de r.
\end{equation}

En sommant (\ref{lemmes:eq1}) et (\ref{lemmes:eq2}), on obtient
\begin{eqnarray}
\lefteqn{\int_0^R\left(\frac{\sinh r}{\cosh r}
a_1'^2+\frac{\cosh r}{\sinh r}a_3'^2\right)\de r}\\
&\geq&\frac{C_1+C_2}{R^2}\inf_{c\in\R}\int_0^R
\left(\frac{\sinh r}{\cosh r}(a_1-c)^2
+\frac{\cosh r}{\sinh r}a_3^2\right)\de r.\nonumber
\end{eqnarray}
L'application du lemme~\ref{lemmes:dod} permet alors de conclure.
\end{demo}

\section{Minoration du spectre à l'aide d'un recouvrement d'ouverts}\label{mc}

 Dans cette partie, les variétés ne seront pas supposées hyperboliques,
ni de dimension~3, et les formes différentielles considérées pourront
être de degré supérieur ou égal à~2. On notera $\mu_{p,i}$ la $i$-ième
valeur propre du laplacien de Hodge en restriction aux $p$-formes
coexactes.
 La proposition~\ref{lemmes:dod} de Dodziuk reste valable 
pour $\mu_{p,i}$ si on impose à $V_i$ d'être un sous-espace des $p+1$-formes.

J.~McGowan a développé dans \cite{mc93} une technique de minoration
du spectre du laplacien de Hodge-de~Rham se basant sur un recouvrement
de la variété par des ouverts sur lesquels le spectre est bien contrôlé.
Cette technique a eu depuis de nombreuses autres applications.

 Cependant, le lemme de McGowan souffre du défaut de ne pas permettre
de minorer la première valeur propre si la cohomologie des intersections
des ouverts du recouvrement est non triviale. Nous allons ici améliorer
ce résultat en montrant qu'on peut toujours minorer la première valeur
propre à condition de contrôler la manière dont interagissent les 
cohomologies des ouverts et de leurs intersections. Le procédé a déjà 
été utilisé dans \cite{ja08} pour une situation très particulière, nous 
allons le généraliser.

 Précisons d'abord les données de la démonstration. On considère un 
recouvrement $(U_i)$ de la variété $M$ par des ouverts dont on supposera 
qu'ils n'ont pas d'intersection d'ordre supérieur ou égal à~3, et on pose
$U_{ij}=U_i\cap U_j$. Sous cette hypothèse, on a une suite exacte de 
Mayer-Vietoris
\begin{equation}\label{mg:mv}
\to H^p(M)\to\bigoplus_i H^p(U_i)\stackrel\delta\to\bigoplus_{i,j}H^p(U_{ij})
\stackrel{d^*}\to H^{p+1}(M)\to,
\end{equation}
cf.~\cite{bt82}, ch.~1, \textsection~1 et~2.

 L'obstruction à minorer la première valeur propre dans la démonstration 
de McGowan est la cohomologie $\oplus H^p(U_{ij})$. 
Plus précisément, à partir d'une forme $\omega\in\Omega^{p+1}(M)$, J.~McGowan
construit un élément $\varphi_{ij}\in\oplus\Omega^p(U_{ij})$ qui est fermé
et montre que s'il est exact, on peut contruire une forme 
$\varphi\in\Omega^p(M)$ telle que $\de\varphi=\omega$ et dont la norme est
contrôlée, puis appliquer la proposition~\ref{lemmes:dod}.

 Par construction, la classe $\varphi_{ij}\in\oplus H^p(U_{ij})$ est
dans l'image de $\delta$ (cette affirmation sera justifiée plus loin). 
L'idée permettant
d'améliorer le résultat de McGowan consiste à utiliser une section de la
flèche~$\delta$, ou une section partielle, définie sur un sous-espace
de l'image de $\delta$. Plus précisément, en identifiant une classe de 
cohomologie avec son représentant harmonique, on se donne un sous-espace 
$E\subset\Ima\delta\subset\oplus\mathcal H^p(U_{ij})$ (où 
$\mathcal H^p(U_{ij})$ désigne l'espace des $p$-formes harmoniques de $U_{ij}$), 
et une application 
$T:E \to\oplus\Omega^p(U_i)$ telle que pour tout $h\in E$, $\de T(h)=0$ et
$\delta [T(h)]=[h]$. On peut alors minorer une valeur propre dont le rang
est alors d'autant plus petit que la dimension de $E$ est grande 
(la première valeur propre si $E=\Ima\delta$). Le lemme qui suit explicite
cette minoration, qui fait intervenir la norme de l'application $T$:

\begin{lem}\label{min:mc}
Soit $(M,g)$ une variété riemannienne compacte de dimension~$n$, $(U_i)$ 
un recouvrement de $M$ par des ouverts n'ayant pas d'intersection d'ordre~3 
ou plus, $p$ un entier $\leq n$ et $\mu(U_i)$ (resp.~$\mu(U_{ij})$) 
la première valeur propre du laplacien restreint aux $p$-formes coexactes 
de $U_i$ (resp. aux $p-1$-formes coexactes de~$U_{ij}$) pour la condition 
de bord absolue. Si on se
donne une partition de l'unité $(\rho_i)$ associée à $(U_i)$, un espace
$E$ et une application $T$ comme ci-dessus, alors
$$\mu_{p,k+1}(M,g)\geq\frac1{\sum_i\left(\frac1{\mu(U_i)}+
\sum_j\left(\frac1{\mu(U_i)}+\frac1{\mu(U_j)}\right)\left(8+
\frac{8c_\rho}{\mu(U_{ij})}+4C_T\right)\right)},$$
avec $k=\dimension\Ima\delta-\dimension E$, $C_T=\|T\|^2$ et 
$c_\rho=\sup_i\|\nabla\rho_i\|_\infty$.
\end{lem}

\begin{demo}
Comme rappelé plus haut, le principe consiste, pour toute forme $\omega
\in\Omega^{p+1}(M)$, à construire une primitive de $\omega$ dont la 
norme est contrôlée et d'appliquer la proposition~\ref{lemmes:dod}.
Cette construction se fait par une chasse dans le complexe de \v Cech-de~Rham
ci-dessous (Cf.~\cite{bt82}):
\begin{equation}\label{demo:diag}
\begin{array}{ccccccccc}
0 & \rightarrow & \Omega^{p+1}(M) & \rightarrow & 
\displaystyle\bigoplus_i\Omega^{p+1}(U_i) & \rightarrow & 
\displaystyle\bigoplus_{ij}\Omega^{p+1}(U_{ij}) & \rightarrow & 0\\
& & \uparrow & & \uparrow & & \uparrow & & \\
0 & \rightarrow & \Omega^p(M) & \rightarrow & 
\displaystyle\bigoplus_i\Omega^p(U_i) & \rightarrow & 
\displaystyle\bigoplus_{ij}\Omega^p(U_{ij}) & \rightarrow & 0\\
& &\uparrow & & \uparrow & & \uparrow & & \\
0 & \rightarrow & \Omega^{p-1}(M) & \rightarrow & 
\displaystyle\bigoplus_i\Omega^{p-1}(U_i) & \rightarrow & 
\displaystyle\bigoplus_{ij}\Omega^{p-1}(U_{ij}) & \rightarrow & 0\\
\end{array}
\end{equation}
 On va d'abord traiter le cas où $E=\Ima\delta$.

Soit $\omega\in\Omega^{p+1}(M)$ une forme exacte. Le but est de construire
une primitive de $\omega$ dont la norme soit contrôlée. On va d'abord
voir comment la construire, puis comment majorer sa norme.

La restriction de $\omega$ 
à chaque $U_i$ est une forme exacte qu'on notera $\omega_i$. 
Si $\varphi_i$ est une primitive de $\omega_i$ pour tout $i$, 
alors les formes $\varphi_{ij}\in\Omega^p(U_{ij})$ définies par
$\varphi_{ij}=\varphi_i-\varphi_j$ sont fermées par définition, et la
construction de la suite de Mayer-Vietoris assure que la classe
$\{[\varphi_{ij}]\}\in\oplus_{ij}H^p(U_{ij})$ est bien dans l'image
de l'application~$\delta$.
 Chaque $\varphi_{ij}$ se décompose donc
en une somme $\varphi_{ij}=\alpha_{ij}+h_{ij}$, où $\alpha_{ij}$ est
exacte et $h_{ij}$ est harmonique (avec condition de bord absolue).
Si pour tout $i$ et $j$, $\beta_{ij}\in\Omega^{p-1}(U_{ij})$ est une 
primitive de $\alpha_{ij}$ (en prenant soin de choisir $\beta_{ji}=
-\beta_{ij}$), on définit $\gamma_i\in\Omega^{p-1}(U_i)$
par $\gamma_i=\sum_j\rho_j\beta_{ij}$. On a alors $\beta_{ij}=\gamma_i-
\gamma_j$, et par commutativité du diagramme (\ref{demo:diag}), on a
$\de\gamma_i-\de\gamma_j=\alpha_{ij}$.

 Par ailleurs, on peut décomposer l'application $T$ en une somme $\sum_iT_i$ 
avec $T_i:\oplus_{ij}H^p(U_{ij})\to\Omega^p(U_i)$. Par conséquent, à chaque
$h_{ij}$, on peut associer la famille $T_i(h_{ij})\in\oplus_i\Omega^p(U_i)$.

Si on pose maintenant $\bar\varphi_i=\varphi_i-\de\gamma_i-\sum_jT_i(h_{ij})$,
on a $\bar\varphi_i-\bar\varphi_j=\varphi_{ij}-\alpha_{ij}-h_{ij}=0$.
Les $\bar\varphi_i$ sont donc les restrictions d'une forme globale
$\bar\varphi\in\Omega^p(M)$ telle que $\de\bar\varphi=\omega$. 
Il reste a majorer la norme de $\bar\varphi$
en fonction de $\omega$ et des données. On commence par écrire
\begin{eqnarray}\label{demo:maj1}
\|\bar\varphi\|^2 & \leq & \sum_i\|\bar\varphi_i\|^2=\sum_i\|\varphi_i\|^2
+\sum_i\|\de\gamma_i+T_i(\sum_j h_{ij})\|^2\nonumber\\
& \leq & \sum_i\|\varphi_i\|^2+2\sum_i\|\de\gamma_i\|^2+
2\sum_i\|T_i(\sum_j h_{ij})\|^2
\end{eqnarray}
en utilisant le fait que $\varphi_i$ est coexacte, donc orthogonale aux 
formes fermées (dans cette inégalité, chaque norme est relative au domaine
sur lequel la forme est définie). Par construction de $\varphi_i$, on a 
\begin{equation}
\|\varphi_i\|^2
\leq\frac{\|\omega_i\|^2}{\mu(U_i)}\leq\frac{\|\omega\|^2}{\mu(U_i)},
\end{equation}
et par ailleurs
\begin{eqnarray}
\|T_i(\sum_j h_{ij})\|^2 & \leq & C_T\sum_j \|h_{ij}\|^2\leq C_T^2\sum_j 
\|\varphi_{ij}\|^2\nonumber\\
& \leq & 2C_T\sum_j(\|\varphi_i\|^2+\|\varphi_j\|^2)\nonumber\\
& \leq & 2C_T\|\omega\|^2\sum_j\left(\frac1{\mu(U_i)}+\frac1{\mu(U_j)}\right).
\end{eqnarray}

 La majoration de $\|\de\gamma_i\|$ s'effectue comme dans \cite{mc93}:
\begin{eqnarray}\label{demo:maj4}
\|\de\gamma_i\|^2 & = & \|\de\sum_j\rho_j\beta_{ij}\|^2\leq
\|\sum_j(\de\rho_j\wedge \beta_{ij}+\rho_j\de\beta_{ij})\|^2\nonumber\\
& \leq & 2\sum_j(c_\rho\|\beta_{ij}\|^2+\|\de\beta_{ij}\|^2)
=2\sum_j(c_\rho\frac{\|\alpha_{ij}\|^2}{\mu_{ij}}+\|\alpha_{ij}\|^2)\nonumber\\
& \leq & 2\sum_j(\frac{c_\rho}{\mu_{ij}}+1)\|\varphi_{ij}\|^2\leq 
4\sum_j(\frac{c_\rho}{\mu_{ij}}+1)
(\|\varphi_i\|^2+\|\varphi_j\|^2).
\end{eqnarray}
En combinant les majorations (\ref{demo:maj1}) à (\ref{demo:maj4}), on 
obtient
\begin{equation}
\|\bar\varphi\|^2\leq\|\omega\|^2\sum_i\left(\frac1{\mu(U_i)}+
\sum_j\left(\frac1{\mu(U_i)}+\frac1{\mu(U_j)}\right)\left(8+
\frac{8c_\rho}{\mu(U_{ij})}+4C_T\right)\right).
\end{equation}
On en déduit une minoration de $\|\omega\|^2/\|\bar\varphi\|^2$ qui donne
le résultat souhaité, en appliquant le lemme~\ref{lemmes:dod}.

Le cas $E\neq\Ima\delta$ se traite comme dans \cite{mc93}: pour tout
espace $V\subset\Omega^{p+1}(M)$ constitué de formes exactes et de 
dimension $k+1$, il existe une forme $\omega\in V$ telle que la forme
harmonique $\sum h_{ij}$ soit dans $E$. On peut alors effectuer les
calculs précédents sur $\omega$ et l'application du lemme~\ref{lemmes:dod} 
donne une minoration de $\mu_{p,k+1}(M,g)$.
\end{demo}

\section{Spectre des variétés hyperboliques}\label{minoration}
On peut maintenant démontrer les deux résultats annoncés dans l'introduction.

\begin{demo}[du théorème~\ref{intro:th1}]
On va appliquer le lemme~\ref{min:mc} à un recouvrement de la variété
par deux ouverts. L'ouvert $U_1$ sera $M_a$ pour un $a<c_{\mathcal M}$ fixé 
(donc $U_1$ englobe la partie épaisse de $M$) et $U_2$ sera la réunion $M_m$ 
des parties minces de $M$. Les valeurs propres $\mu(U_1)$ et $\mu(U_{12})$ 
sont uniformément minorés d'après le lemme~\ref{lemmes:1} (les mêmes 
arguments pemettent de majorer uniformément la constante $c_\rho$). 
Celui de $\mu(U_2)$ étant minoré grâce au lemme~\ref{lemmes:2}, on est 
essentiellement ramené au contrôle de la constante $C_T$ dans le 
lemme~\ref{min:mc}.

On va commencer par la minoration de la $(k+1)$-ième valeur propre. On
définit $E$ comme étant l'image de l'application 
$H^1(U_1)\to H^1(U_{12})$, et $T$ comme étant la réciproque de cette
application, en identifiant chaque classe de cohomlogie à son représentant
harmonique.
Comme la géométrie de $U_1$ est contrôlée, la norme $C_T$ est uniformément
majorée, et d'après le lemme~\ref{topo:image} la dimension de $E$ est $k$.
Le lemme~\ref{min:mc} donne donc la minoration
\begin{equation}\label{minoration:eq2}
\mu_{k+1}(M,g)\geq\frac C{C'+R^2}
\end{equation}
où $C,C'>0$ sont deux constantes ne dépendant que de $V$ et $R^2$ étant
le plus grand rayon des composantes connexes de $M_m$. Comme le diamètre
$d$ de la variété est uniformément minoré (on a nécessairement $d\geq
c_{\mathcal M}$) et que $d>R$, on peut trouver une constante $c$ telle que 
$\mu_{k+1}(M,g)\geq\frac c{d^2}$. 

 Pour la minoration de la première valeur propre, nous allons d'abord 
traiter le cas où la partie mince de la variété ne compte qu'une seule 
composante connexe, c'est-à-dire que $k=1$. On verra ensuite comment la 
démonstration se généralise.

 Si $k=1$, l'espace de cohomologie 
$H^1(U_{ij})=H^1(U_{12})$ est alors de dimension~2, et les images des deux 
applications $H^1(U_i)\to H^1(U_{12})$ sont de dimension~1.
\begin{figure}[h]
\begin{center}
\begin{picture}(0,0)%
\includegraphics{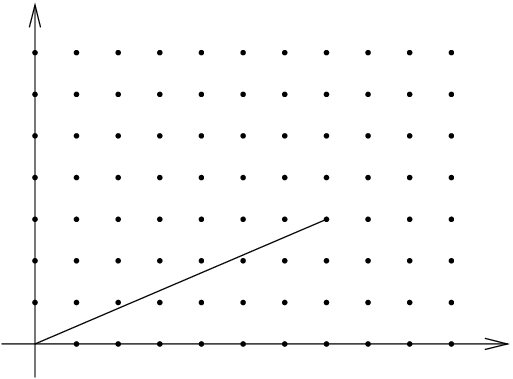}%
\end{picture}%
\setlength{\unitlength}{4144sp}%
\begingroup\makeatletter\ifx\SetFigFont\undefined%
\gdef\SetFigFont#1#2#3#4#5{%
  \reset@font\fontsize{#1}{#2pt}%
  \fontfamily{#3}\fontseries{#4}\fontshape{#5}%
  \selectfont}%
\fi\endgroup%
\begin{picture}(3894,2878)(278,-2956)
\put(2623,-1672){\makebox(0,0)[lb]{\smash{{\SetFigFont{12}{14.4}{\rmdefault}{\mddefault}{\updefault}{\color[rgb]{0,0,0}$(a,b)$}%
}}}}
\end{picture}%
\end{center}
\caption{\label{minoration:diag}}
\end{figure}

La figure~\ref{minoration:diag} représente $H^1(U_{12})$, en choisissant
l'image de $H^1(U_1)$ comme axe des abscisses. Comme les deux applications
$H^1(U_i)\to H^1(U_{12})$ envoient classe entière sur classe entière,
l'image de $H^1(U_2)$ est déterminée par un couple d'entier $(a,b)$
qui est l'image d'un générateur $g_2$ de $H^1(U_2,\Z)$ (par exemple la
forme $\frac1l\de t$ de l'expression~(\ref{topo:g})). Le couple
$(a,b)$ correspond, à un changement de base près, au paramètre
$(p,q)$ défini en~(\ref{topo:chir}) et caractérisant la topopogie de la 
variété à partie épaisse fixée. Il est
possible que $b=0$, mais comme il n'y a qu'une seule variété
(à partie épaisse fixée) vérifiant cette propriété, on peut exclure
ce cas sans nuire à la généralité. Si de plus 
on note $g_1\in H^1(U_1)$ un élément dont l'image est $(1,0)$, on
pose $E=H^1(U_{12})$ et on cherche une section $T:E\to\langle g_1,g_2\rangle$.

La matrice de $\oplus H^1(U_i)\to H^1(U_{12})$ restreinte à 
$\langle g_1,g_2\rangle$ est $\smat{1&a\\0&b}$. Celle de $T$ est son
inverse 
\begin{equation}
T :\left(\begin{array}{cr}
1 & -\frac ab \\ 0 & \frac1b
\end{array}\right).
\end{equation}
On veut majorer la norme de $T$ en fonction de celle de sa matrice. Pour 
cela, on doit contrôler la norme de chaque composante $T_i: H^1(U_{12})\to 
 H^1(U_i)$ de $T$ en fonction de celle de la colonne de $T$ 
correspondante. Comme les métriques sur $U_1=M_a$ et $U_{12}$ sont contrôlées,
la norme de $T_1: H^1(U_{12})\to H^1(U_1)$ est uniformément majorée. Dans 
le cas de $U_2$, on doit
estimer la norme de la forme harmonique $\de t$ pour la métrique
(\ref{lemmes:g}) en fonction de la norme de cette forme en restriction
à $U_{12}$. Or
\begin{equation}
\|\de t\|_{U_2}^2=\int_{U_2} |\de t|^2\de v_g=\int_{U_2}\cosh^{-2}r\de v_g
=2\pi l\int_0^R\tanh r\de r
\end{equation}
et 
\begin{equation}
\|\de t\|_{U_{12}}^2=\int_{U_{12}} |\de t|^2\de v_g
=2\pi l\int_{R_a}^R\tanh r\de r,
\end{equation}
où $R_a$ est la coordonnée du bord de $M_a$ dans la partie mince. 
$R_a$ dépend de $|R-R_a|$ mais converge quand $R\to+\infty$, donc il existe
une constante $c>0$ telle que $\|\de t\|_{U_2}^2/\|\de t\|_{U_{12}}^2\leq cR$. 
On obtient donc qu'il existe une constante $c'>0$ telle que 
\begin{equation}\label{minoration:eq1}
C_T\leq c'\max(1,R|\frac ab|^2)\leq c'R a^2.
\end{equation}

Pour finir, on remarque que $a^2+b^2$ est asymptotiquement de l'ordre
de $|\de t/l|^2$, cette norme étant calculée le long du bord du tube 
(c'est-à-dire qu'on prend $r=R$ dans l'expression (\ref{topo:g}) de 
la métrique), et donc que $a^2+b^2\sim l^{-2}e^{-2R}$. On en 
déduit, en utilisant la relation (\ref{topo:circ}), qu'il existe une constante
$c''$ telle que $a^2\leq c''e^{2R}$, et donc $C_T\leq c'c''Re^{2R}$. 
En appliquant le lemme~\ref{min:mc}, on obtient l'existence de constantes
$C_i$, $i=1,2,3$ telles que
\begin{equation}
\mu_1(M)\geq\frac1{(C_1+R^2)(C_2+C_3Re^{2R})}
\end{equation}
Le théorème s'en déduit.

 On passe maintenant au cas d'un nombre $k$ quelconque de parties minces.
L'espace $H^1(U_{12})$ est alors de dimension~$2k$, et les images des deux
applications $H^1(U_i)\to H^1(U_{12})$ sont de dimension~k. On ne peut 
pas exclure le cas où ces deux images ne sont pas en somme directe 
(cas $b=0$ dans ce qui précède). On va donc noter $E$ l'image de 
$H^1(U_1)\oplus H^1(U_2)\to H^1(U_{12})$ et $l$ sa dimension.
 On note en outre $g_1,\ldots,g_k$ une base du supplémentaire $F$ de 
$\Ker(H^1(U_1)\to H^1(U_{12}))$ formée de classes entières, et 
$g_{k+1},\ldots,g_{2k}$ des générateurs de chacune des composantes 
de la partie mince.

On va maintenant écrire la matrice de $F\oplus H^1(U_2)\to E$. La matrice
$A$ de $F\to E$ a $k$ colonnes de $l$ lignes, et l'image de 
chaque $g_{k+i}$ est un vecteur à $l$ composantes qu'on notera $B_i$. 
On obtient une matrice de la forme
\begin{equation}\label{minoration:mat1}
\left(\begin{array}{c|cccc}
A & B_1 & B_2 & \ldots & B_k
\end{array}\right).
\end{equation}
Cette matrice à $2k$ colonnes et $l$ lignes, elle n'est donc pas 
nécessairement carrée. Mais comme l'application $F\to E$ est injective, on
peut rendre la matrice~(\ref{minoration:mat1}) carrée et inversible en
enlevant $2k-l$ des colonnes $B_i$. On obtient la matrice
\begin{equation}\label{minoration:mat2}
P=\left(\begin{array}{c|cccc}
A & B'_1 & B'_2 & \ldots & B'_{l-k}
\end{array}\right),
\end{equation}
les vecteurs colonnes $B'_1,\ldots B'_{l-k}$ étant l'image d'une sous familles 
$g'_1,\ldots g'_{l-k}$ de $g_{k+1},\ldots g_{2k}$.

On choisit l'application $T$ comme étant l'inverse de celle définie par la
matrice~$P$. Elle va donc de $E$ dans $F\oplus\langle
g'_1,\ldots g'_{l-k}\rangle$ et sa matrice est $P^{-1}$. Comme $P$ est
à coefficients entiers, que $|\det P|\geq1$ et que 
$P^{-1}=(\det P)^{-1}\cdot \overline P$ où
$\overline P$ désigne la comatrice de $P$, il suffit de majorer la  
norme de $(\det P)\cdot T$ dont la matrice est $\overline P$.
Or, les coefficients de $\overline P$ sont des polynômes de degré~1 par
rapport à chacun des $B'_i$, on peut donc écrire
\begin{equation}
\|\overline P\|=c\cdot\prod_{i=1}^{k-l}\|B'_i\|,
\end{equation}
où $c$ est une constante ne dépendant que de $A$.

Comme dans le cas $k=1$, la norme de la composante $T_1:E\to H^1(U_1)$
est majorée par celle de $\overline P$, et chacune des composantes
$T'_i:E\to\langle g'_i\rangle$ est majorée par $c'R\|\overline P\|$ où 
$c'$ ne dépend que de la topologie de la partie épaisse.
Pour finir, toujours comme dans le cas $k=1$, la norme $\|B'_i\|$ est 
majorée par $c''e^{R_i}$, où $R_i$ est le rayon de la partie mince 
correspondante.

On obtient finalement que 
\begin{equation}
C_T\leq CR\|\overline P\|\leq C'Re^{2(k-l)R},
\end{equation}
où $C$ et $C'$ ne dépendent que de la topologie de la partie épaisse.
On conclut ensuite comme dans le cas $k=1$.
\end{demo}

\begin{demo}[du théorème~\ref{intro:th2}]
Si la variété $M$ a plusieurs cusps, on commence par tous les remplir sauf
un en utilisant le théorème de Thurston (\cite{gr81}, théorème~4.A) et
on note $M'$ la variété obtenue. Un remplissage du dernier cusp est
alors paramétré par un élément $(p,q)\in\mathcal P$ auquel correspond
un élément $(a,b)$ de $H^1(T^2)$ comme sur la figure~\ref{minoration:diag}.

On peut alors appliquer les résultats de la démonstration précédente à
la variété compacte $M''$ ainsi obtenue en remplaçant $M''_{c_{\mathcal M}}$
(resp. $M''_a$), par le même domaine de $M''$ auquel on ajoute les 
composantes de la parties minces obtenues par le remplissage des premiers
cusps. En particulier, si on choisit $|b|\geq R|a|$, 
l'inégalité~(\ref{minoration:eq1}) donne une majoration uniforme de $C_T$,
auquel cas la minoration~(\ref{minoration:eq2}) s'applique à la première
valeur propre.

Il suffit donc de choisir une suite $M_i$ de variétés données par 
la suite $(a_i,b_i)=(1,i)$. Pour $i$ assez grand la variété $M_i$
est bien hyperbolique, on a $a_i^2+b_i^2\sim e^{2R}$ donc $b_i\sim
e^R$ et la condition  $|b|\geq R|a|$ est bien vérifiée.
\end{demo}

\noindent Pierre \textsc{Jammes}\\
Université d'Avignon et des pays de Vaucluse\\
Laboratoire d'analyse non linéaire et géométrie (EA 2151)\\
F-84018 Avignon\\
\texttt{Pierre.Jammes@ens-lyon.org}
\end{document}